\documentclass[10pt]{article}
\usepackage{amsfonts}
\usepackage{amsfonts}  
\usepackage{pdfsync}
\title{Radial singular solutions of fully nonlinear equations in punctured balls}

\author{Isabeau Birindelli \thanks{isabeau@mat.uniroma1.it}\\
Dipartimento di Matematica, Sapienza Universit\`a\  di Roma
\and
  Fran\c{c}oise Demengel\thanks{francoise.demengel@cyu.fr}\\
  D\'epartement de Math\'ematiques,
CY-Cergy Paris Universit\'e\  
  \and
   Fabiana  Leoni\thanks{leoni@mat.uniroma1.it}\\ 
   Dipartimento di Matematica, Sapienza Universit\`a\  di Roma}

\date{}

\usepackage{amsthm}
\usepackage{amsmath}
\usepackage{color}

\catcode`@=11 \@addtoreset{equation}{section} \catcode`@=12

\newtheorem{theo}{Theorem}[section]
\newtheorem{prop}[theo]{Proposition}
\newtheorem{rema}[theo]{Remark}

\newcommand{\dis}{\displaystyle}

\def\R{\mathbb  R}

\begin{document}

\maketitle
\begin{abstract}
\noindent
We study fully nonlinear uniformly elliptic equations having a singular reaction term with inverse quadratic potential and an absorbing superlinear  term of $p$-power type. We consider  equations posed in  punctured balls centered at the origin,  and we prove that all radial solutions are singular around the origin, by  providing a complete classification in dependence of $p$ of  their   asymptotic behavior  near the singularity.
\end{abstract}

\vspace{1cm}

\noindent
\emph{2020 Mathematical Subject Classification }: 35J60, 35J75, 35B40.

\medskip
\noindent
\emph{Keywords}:  Fully nonlinear elliptic equations; singular potentials; superlinear absorbing terms.

\bigskip 

 Research partially supported by GNAMPA-INdAM Project 2023 ”Equazioni completamente non lineari locali e non locali” CUP E53C22001930001 and PRIN 
2022 ”PDEs and optimal control methods in mean field games, population dynamics and multi-agent models”

 \section{Introduction }
We are interested in the existence and the behavior around the origin of radial  positive solutions of fully nonlinear elliptic equations posed in punctured balls centered at the origin, namely equations of the form 
\begin{equation}\label{eqmain1}
F( D^2 u) + \mu \frac{u}{|x|^2} = u^p\, .
\end{equation}
$F$ will be one of the Pucci's extremal operators, although some of the results can be proved for more general fully nonlinear elliptic operators.  The exponent $p$ will be  supposed to be larger than one, that is we put ourselves in a superlinear regime, and the parameter $\mu$ will be  assumed to be positive and smaller than the principal eigenvalue associated with Pucci's operators and the inverse quadratic potential $r^{-2}$, as recalled below.
\medskip

For $0< \lambda\leq  \Lambda$ given,  Pucci's sup-operator is defined as
\begin{equation}\label{M+}
 \mathcal {M}^+ ( M) =\mathcal{M}^+_{\lambda, \Lambda} (M)= \Lambda \sum_{\mu_i>0}\mu_i+\lambda\, \sum_{\mu_i<0}\mu_i\, ,
 \end{equation}
where $\mu_1,\ldots,\mu_N$ are the eigenvalues of any  matrix $M$ belonging to the space $\mathcal{S}_N$ of the $N-$squared symmetric matrices. Equivalently, $ \mathcal {M}^+$ can be defined as
$$ \mathcal {M}^+ ( M)=  \sup_{A\in {\mathcal A}_{\lambda,\Lambda}} {\rm tr} (AM)$$
where ${\mathcal A}_{\lambda,\Lambda}=\left\{A\in{\mathcal S}_N\,:\,\lambda \, I_N\leq A\leq \Lambda \, I_N\right\}$, $I_N$ being the identity matrix.

Symmetrically,  Pucci's inf-operator is defined as
 \begin{equation}\label{M-} 
 \mathcal {M}^- ( M) =\mathcal{M}^-_{\lambda, \Lambda} (M)= \lambda \sum_{\mu_i>0}\mu_i+\Lambda\, \sum_{\mu_i<0}\mu_i
=\inf_{A\in {\mathcal A}_{\lambda,\Lambda}} {\rm tr} (AM)\, ,
\end{equation}
and  they are mutually related by 
$$
\mathcal {M}^+ ( M)=-\mathcal {M}^- (- M)\qquad \forall\, M\in \mathcal{S}_N\, .
$$
Pucci's extremal operators have been introduced in 
\cite{Pucci}, and then extensively studied in \cite{CC}. Acting as barriers not only for linear operators, but in the 
 whole class of uniformly elliptic operators with the same ellipticity constants, they play a crucial role in the regularity theory for fully nonlinear uniformly elliptic equations, see \cite{CC}. Moreover, as sup and inf of linear operators, they  frequently arise in 
 the context of  optimal stochastic control problems, see e.g. \cite{FS, L}, with special application to mathematical finance.
		
Pucci's extremal operators  can be seen as a generalization of  Laplace operator, since
$${\mathcal M}^-_{\lambda,\lambda}(X) ={\mathcal M}^+_{\lambda,\lambda}(X)=\lambda\, {\rm tr}(X)\, .$$
When $F$ is the Laplace operator, problem \eqref{eqmain1} has been extensively studied, also in the case $\mu=0$, see \cite{Ci, CirsChau, CDu1, CDu2, WDu}.  However, the employed techniques rely on the variational structure of the problem, and they cannot be immediately extended to the fully nonlinear framework.

As soon as  $\lambda<\Lambda$,  the operators   ${\mathcal M}^\pm_{\lambda,\lambda}$ loose the   linear structure and the   divergence form. However, nice structure properties can be partially recovered when Pucci's operators are evaluated on radially symmetric functions. Indeed, if $u(x)=u(|x|)=u(r)$ is a radial function of class $C^2$ for $x\neq 0$, then the eigenvalues of the hessian matrix $D^2u(x)$ are given by $u''(r)$, which is simple, and $\frac{u'(r)}{r}$, which has multiplicity $N-1$. Thus, if we know the sign of $u'(r)$ and $u''(r)$, we can explicitly compute the expression of $\mathcal{M}^\pm(D^2u)$ according to definitions \eqref{M+} and \eqref{M-}. In particular, if $u(x)=u(|x|)$ is convex and decreasing as a function of $r=|x|$, then one has
$$
\mathcal{M}^+(D^2u)= \Lambda \, \left( u''+\frac{\lambda}{\Lambda} (N-1) \frac{u'}{r}\right)
$$
as well as
$$
\mathcal{M}^-(D^2u)= \lambda \, \left( u''+\frac{\Lambda}{\lambda} (N-1) \frac{u'}{r}\right)\, .
$$
These expressions motivate the introduction of the dimension like parameters
 $$
\tilde{N}_+= \frac{\lambda}{\Lambda} (N-1)+1\, ,\qquad \tilde{N}_-= \frac{\Lambda}{\lambda} (N-1)+1
$$
associated with $\mathcal{M}^+$ and $\mathcal{M}^-$ respectively. With such a position, we obtain the expressions
$$
\begin{array}{c}
\mathcal{M}^+(D^2u)= \Lambda \, \left( u''+(\tilde{N}_+-1)\displaystyle  \frac{u'}{r}\right)\\[3ex]
\mathcal{M}^-(D^2u)= \lambda \, \left( u''+(\tilde{N}_--1) \displaystyle \frac{u'}{r}\right)
\end{array}
$$
which show that, when evaluated on smooth radially symmetric functions $u$, Pucci's operators take an expression structurally identical to that of Laplace operator, provided that the monotonicity and convexity properties of $u$ are known. More in general, radially symmetric solutions of equations involving Pucci's operators satisfy associated ordinary differential equations having piecewise constant coefficients jumping at the points where solutions change monotonicity and/or convexity.
\medskip

The purely eigenvalue problem associated with Pucci's operators and the singular potential $r^{-2}$, i.e. the Dirichlet boundary value problem 
$$
\left\{
\begin{array}{c}
\mathcal{M}^\pm (D^2u) + {\bar{\lambda}}^\pm \displaystyle \frac{u}{r^2}=0 \quad \hbox{ in } B_1\setminus \{0\}\\[2ex]
u=0 	\quad \hbox{ on } \partial B_1	
\end{array}
\right.
$$
where $B_1$ is the unit ball centered at zero, has been recently studied in \cite{BDLarxiv}. We explicitly 	determined the principal eigenvalues ${\bar \lambda}^\pm$, that are the only values for which the above problems admit positive solutions, and the related radial principal eigenfunctions, which are shown to be unique up to multiplicative constants. Precisely,  one has 
$$
\bar{\lambda}^+ =\Lambda \left(\frac{\tilde N_+-2}{2}\right)^2\, ,\quad  \bar{\lambda}^- =\lambda \left(\frac{\tilde N_--2}{2}\right)^2\, ,
$$
and the related eigenfunctions are given by
$$
\Phi^\pm(r)=  \frac{-\log r}{r^{\frac{\tilde{N}_\pm-2}{2}}}\, .
$$
 In the sequel we  will focus on the operator  $\mathcal{M}^+$ only, the changes to be made for $\mathcal{M}^-$ being obvious.
 
For notational simplicity, we will denote merely by $\bar \lambda$ the principal eigenvalue $\bar{\lambda}^+$, and we will set
$$
\tau\, := \frac{\tilde{N}_+ -2}{2}\, ,
$$
so that
$$
\bar \lambda= \Lambda \tau^2\, .
$$
In the present paper, we will be concerned with the existence and the asymptotic behavior near the origin of positive radial solutions of the equation
  \begin{equation}\label{eqmain}
  {\mathcal M}^+( D^2 u) + \mu \frac{u}{r^2} = u^p\, .
 \end{equation}
 By exploiting the invariance of the equation with respect to the scaling
 $$
 u_\alpha (r)=\alpha^{\frac{2}{1-p}} u \left( \frac{r}{\alpha}\right)\qquad \hbox{ for all } \alpha>0\, ,
 $$
 we will consider without loss of generality equation \eqref{eqmain} posed in 
 $ B_1\setminus \{0\}$, 
 and we will refer to smooth, positive, radial solutions  simply as solutions. 
  
  In order to state our main results, let us introduce the two constants
  $$
  \tau^\pm = \tau \pm \sqrt{\tau^2 -\frac{\mu}{\Lambda}}\, ,
  $$
  which are well defined under the assumption $\mu \leq \bar{\lambda}$ and correspond to the roots of the algebraic equation
  $$
  x^2 +\left( \tilde{N}_+-2\right) x +\frac{\mu}{\Lambda}=0\, .
  $$
The existence and the asymptotic behavior of solutions of equation \eqref{eqmain} clearly depend on the mutual values   of the parameter $\mu$ and $p$ with respect to $\bar \lambda$ and $\tau^\pm$ respectively. In particular, let us  introduce the critical values
\begin{equation}\label{pstar}
p^\star = 1+ {2 \over \tau^+}
\end{equation}
and 
 \begin{equation}\label{pstarstar}
p^{\star \star}= 1+ {2 \over \tau^-}\, ,
\end{equation}
which satisfy, for $0<\mu < \bar \lambda$, 
$$
1<p^\star <p^{\star \star}\, .
$$
Our first main result reads as follows.

     \begin{theo}\label{mainresult}
Let  us assume $0<\mu <\bar{\lambda}$ and $p>1$, and let  $p^\star$ and  $p^{ \star \star}$ be  defined as in \eqref{pstar} and \eqref{pstarstar} respectively. Then:    
\begin{enumerate}
     \item 
           if $p < p^{ \star }$,  then equation \eqref{eqmain}  has at least a  solution $u$ satisfying 
  \begin{equation}\label{2surp-1} 
   r^{2\over p-1}   u(r) \to  K \qquad \hbox{ as } r\to0\, ,
        \end{equation}
 at least a  solution $v$ satisfying
 \begin{equation}\label{tau+}
    r^{\tau^+}  v(r) \to  c_1 \qquad \hbox{ as } r\to 0\, ,
\end{equation}
and at least a solution $w$ satisfying
  \begin{equation}\label{tau-}
 r^{\tau^-}  w(r) \to  c_2  \qquad \hbox{ as } r\to0\, ,
        \end{equation}        
where   $K \, := \left[ \Lambda \left( {2\over p-1}-\tau^+\right) \left( {2\over p-1}-\tau^-\right)\right]^{1\over p-1}$ and  $c_1,\ c_2>0$ are suitable  constants;
      
   \item if  $p^\star\leq p< p^{ \star \star}$, then equation \eqref{eqmain} has at least a   solution $u$ satisfying
 \eqref{tau-};
 
 \item if $p = p^{ \star \star}$, then  equation \eqref{eqmain}  has at least a   solution $u$ satisfying 
\begin{equation}\label{tau-log}
 r^{\tau^-} ( -\log r)^{\tau^-/2}u(r) \to \bar K \qquad \hbox{ as } r\to 0\, ,
 \end{equation}
        where 
$\bar K = \left[\Lambda \tau^- (\tau -\tau^-)\right]^{\tau^-/2}$;
        
        \item if  $p> p^{\star\star} $,  then  equation \eqref{eqmain}  has at least a    solution $u$ satisfying (\ref{2surp-1}).
 \end{enumerate}
 \end{theo}

In the next result we  show that the asymptotic behaviors around the origin detected in Theorem \ref{mainresult}  are actually the only possible ones for any solution.

 \begin{theo}\label{behave}
Let  $u$ be a solution of equation \eqref{eqmain} with $p>1$ and $0<\mu <\bar \lambda$. Then:
\begin{enumerate}

\item if $1<p < p^{ \star }$, then  $u$ satisfies either \eqref{2surp-1} or  \eqref{tau+} or 
 \eqref{tau-};

\item if  $p^\star \leq p < p^{ \star \star}$,  then    $u$ satisfies    \eqref{tau-};  

\item if $p = p^{ \star\star}$, then $u$ satisfies \eqref{tau-log};

\item  if $p>p^{\star \star}$, then $u$ satisfies \eqref{2surp-1}.
 \end{enumerate}            
             \end{theo}
Theorems \ref{mainresult} and \ref{behave} will be obtained as consequences of some preliminary results and the combination of several techniques.   

In Section 2 we give two comparison theorems. The first one holds for not necessarily radial solutions of the general equation \eqref{eqmain1} posed in a domain $\Omega$ far away from the singular point $0$, and simply follows from the ellipticity of the principal part and the superlinear character of the absorbing term. In the second result, we remarkably extend the comparison property to sub- and super-solutions of equation \eqref{eqmain} posed in a ball centered at the singular point, under a growth condition on the functions to be compared. Although limited to radial solutions, this second  comparison result shows how, for any $p>1$,   the superlinear absorbing term can be used  in order to   balance the singular reaction one, see next Theorem \ref{compa2}, and it will play a crucial role in the proof of Theorem \ref{behave}. 
\smallskip

Theorem \ref{mainresult} will be proved in Section 3. We first construct explicit sub and super-solutions to be used as barrier functions. Then, the  different solutions with the desired asymptotic behavior near the origin will be obtained as a result of an iterative procedure applied on truncated annular domains,  combined with the comparison property on domains not containing the singular point.
 \smallskip
             
The final Section 4 will be devoted to the proof of Theorem \ref{behave}. 

For the proof of statements 1 and 4 we will make use of the so called Emden--Fowler transform, i.e.  we perform the change of unknown 
$$
x(t)=r^{\frac{2}{p-1}}u(r)\, ,\quad r=e^t\, ,
$$
which reduces the initial problem \eqref{eqmain} for $u(r)$ to an autonomous equation for $x(t)$, see equation \eqref{eqx}. Initially introduced in \cite{F} for the study of semilinear elliptic equations with reaction terms of  power-type,  the Emden-Fowler transform has been  applied in the analysis of radial solutions of  fully nonlinear equations for the first time in \cite{FQ}, and subsequently used in \cite{BGLP, GIL, GILP}. It fits perfectly also to study equation \eqref{eqmain}, owing to its scaling property. The different asymptotic behaviors \eqref{2surp-1}, \eqref{tau+} and \eqref{tau-} will be derived as consequences of the convergence of $x(t)$ to the different equilibria of the equation it satisfies. In particular, the convergences \eqref{tau+} and \eqref{tau-} for $u$ correspond to the convergence  
$$x(t)\to 0\qquad \hbox{ as } t\to -\infty\, ,
$$
with $x(t)$ being  asymptotic to either one of the solutions of the linearized equation associated with the equation satisfied by $x$.

Even if  statement 2 also could be proved by using the properties of the Emden-Fowler transform,  we will provide an alternative proof based on  a different change of variable. In this case, we introduce the new function
$$
v(r)=r^{\tau-} u(r)\, ,
$$             
and we prove that $v(r)$ is monotone non decreasing for $r>0$ and it satisfies an equation which enjoys the comparison principle. 

The function $v$ will be used also in the delicate proof of statement 3, which  will be finally deduced by constructing  sharp sub- and super-solutions of the equation satisfied by $v$.
\medskip

In conclusion, let us finally mention that the results of the present paper can be regarded as a first step into the study   of fully nonlinear uniformly elliptic equations having   coefficients with isolated singularities. Many  questions remain open, ranging from the analysis of non radially symmetric solutions to the extension of the results to more general elliptic operators and zero order terms. In particular, the study of the limit cases $\mu=0$ and  $\mu=\bar \lambda$ is  postponed to a future work.

        \section{Comparison Principles}   
             
  In this Section we provide   two  comparison theorems: the former is a classical comparison result, for equations having  continuous coefficient up to the boundary, which will be applied in annular domains; the latter is a global result, which gives a comparison result up to the singularity point, under a growth assumption on the sub and super solutions to be compared. 
             
 \begin{prop}\label{compa1}
Let $p>1$, and 
 suppose that $u$ and $v$ are positive functions satisfying in the classical sense
 $$ \mathcal{M}^+( D^2 u) + \mu \frac{u}{|x| ^2} - u^p \geq 0 \geq 
 \mathcal{M}^+( D^2 v)+ \mu \frac{v}{|x| ^2} -v^p $$
  in  a domain $\Omega$ such that $0\notin \overline{\Omega}$. We suppose also that 
  $u$ and $v$ are positive and continuous in  $\overline{\Omega}$, with $u\leq v$ on $\partial \Omega$. Then, $u\leq v$ in $\overline{\Omega}$. 
  \end{prop}
  
\begin{proof} 
    Since $v>0$ in $\overline{ \Omega}$, then the supremum
    $$\gamma = \sup_{\overline \Omega} { u \over v}$$
     is well defined and it is indeed a maximum. Let $\bar x\in \overline{\Omega}$ be a maximum point. If $\bar x\in \Omega$, then  one has
$$ \gamma^p v^p( \bar x)-\mu \gamma \frac{v( \bar x)}{|\bar x|^{2}} = u^p(\bar x)- \mu \frac{u(\bar x)}{|\bar x|^{2}} \leq \mathcal{M}^+( D^2 u(\bar x)) \leq \mathcal{M}^+ \left( D^2  (\gamma v ( \bar x))\right)\leq \gamma v^p( \bar x)-\mu\gamma \frac{v( \bar x)}{ |\bar x|^{2}}$$
which  yields  $\gamma \leq1$. On the other hand, if $\bar x\in \partial \Omega$, then $\gamma \leq 1$ as well by assumption. Thus, in any case, we get $\gamma\leq 1$.
      
       \end{proof} 
       
 \begin{rema} 
{\rm    We can replace the zero order nonlinearity $u\mapsto u^p$ by $u\mapsto g(u)$ for any continuous function $g$ such that ${g(u)\over u}$ is strictly increasing. Moreover, the coefficient   $\frac{\mu}{|x|^2}$ can be replaced by any continuous and bounded   function, and the proof can be extended to merely semicontinuous functions satisfying the differential inequalities in the viscosity sense.}
   
 \end{rema}  
 
\medskip
   
\begin{theo}\label{compa2}
Let  $\mu>0$, $p>1$, $r_0>0$ and let $u$ and $v$ be two  positive radial  functions in $C^2\left(B_{r_0}\setminus \{0\}\right)$ satisfying
$$
\mathcal{M}^+(D^2u)+\mu \frac{u}{r^2} -u^p\geq 0 \geq \mathcal{M}^+(D^2v) +\mu \frac{v}{r^2} -v^p\quad \hbox{ in } B_{r_0}\setminus \{0\}\, .
$$
Assume that $u\, ,\ v \in C\left( \overline{B_{r_0}}\setminus \{0\}\right)$ and that
 there exist positive constants $c_1, c_2$ such that 
 $$ c_1 r^{-2\over p-1}\leq  v(r)\quad \hbox{ and }\quad   u(r)\, , v(r) \leq c_2 r^{-2\over p-1}\quad \hbox{  for }\  0<r\leq r_0\, .$$
 If $u(r_0)\leq v(r_0)$, then $u\leq v$ in $\overline{B_{r_0}}\setminus \{0\}$.
             \end{theo}
             
              \begin{proof} 
Let us consider the supremum
$$
\gamma'= \sup_{0<r\leq r_0} \frac{u(r)}{v(r)}\, ,
$$
which is well defined by the growth assumptions on $u$ and $v$. We assume, by contradiction, that $\gamma' >1$ and we select a maximizing sequence $\{ r_n\} \subset  (0,r_0]$ converging to a point $0\leq \bar r\leq r_0$. If  $\bar r>0$, we easily reach the contradiction $\gamma'\leq 1$ by arguing as in the proof of Proposition \ref{compa1}. Thus, let us consider the case $r_n\to 0$ as $n\to \infty$. Without loss of generality, we can assume that $\{r_n\}$ is decreasing and that, for any $n\geq 1$, $u(r_n)> \gamma \, v(r_n)$ for some $1<\gamma<\gamma'$ arbitrarily fixed.

 We claim that   there exists $n_1\geq 1$ such that $u(r)> \gamma v(r)$ for all $0<r\leq r_{n_1}$. Indeed, if not, we could find, by continuity,  some interval $[\sigma, s] \subset (0, r_1]$ such that $u(r)-\gamma v(r)>0$ for $\sigma<r<s$ and $u(\sigma)-\gamma v(\sigma)= u(s)-\gamma v(s)=0$. On the other hand, in the interval $(\sigma, s)$ we have
 $$
 \begin{array}{rl}
 \dis \mathcal{M}^+( D^2 ( u-\gamma v) ) & \dis \geq \mathcal{M}^+(D^2u) -\gamma \mathcal{M}^+(D^2v)\\[2ex]
  & \dis \geq u^p-\gamma  v^p -\mu \frac{u-\gamma v}{r^2} \\[2ex]
  & \dis \geq ( \gamma^p-\gamma) v^p -\mu (\gamma' -\gamma)\frac{v}{r^2}  \geq \left[ (\gamma^p-\gamma) c_1^p -\mu (\gamma'-\gamma)c_2\right] r^{\frac{-2p}{p-1}}\, .
  \end{array}$$         
By choosing $\gamma>1$ sufficiently close to $\gamma'$,  the right hand side of the above inequality can be made strictly positive. On the other hand, the function $u-\gamma v$ is positive in $(\sigma, s)$ and zero on the boundary, hence it has an interior maximum point $t$ at which $(u-\gamma v)'(t)=0$ and $(u-\gamma v)''(t)\leq 0$, so that $\mathcal{M}^+\left( D^2 (u-\gamma v)(t)\right)\leq 0$.

 The above argument  shows not only that $u(r)- \gamma v(r)>0$ for all $0<r\leq r_{n_1}$, but also that  
 $$
\mathcal{M}^+( D^2 ( u-\gamma v) )>\frac{c_3}{r^{\frac{2p}{p-1}}} \quad \hbox{ for } 0< r<r_{n_1}\, ,
$$ 
with e.g. $c_3=\frac{(\gamma^p-\gamma) c_1^p}{2}>0$. 
Moreover, since $u-\gamma v$ cannot have any local maximum point in $(0,r_{n_1})$, it follows  that $(u-\gamma v)'$ has constant sign in a sufficiently small right neighborhood of zero (see also Lemma 2.1 in \cite{BDLarxiv}). By observing that
 $$
 \lim_{n\to \infty} \left( u(r_n)-\gamma v(r_n)\right) = \lim_{n\to \infty} v(r_n) \left( \frac{u(r_n)}{v(r_n)}-\gamma\right)= +\infty\, ,
 $$
 we deduce that $(u-\gamma v)'<0$ for $r$ sufficiently small. Hence, one has
$$
\frac{c_3}{r^{\frac{2p}{p-1}}}<\mathcal{M}^+(D^2(u-\gamma v))=\Lambda (u-\gamma v)''+ \lambda (N-1) \frac{(u-\gamma v)'}{r} \, ,
 $$
that is 
$$
\left( r^{\tilde{N}_+-1}(u-\gamma v)'(r) \right)'>\frac{c_3}{\Lambda\, r^{\frac{2p}{p-1}-\tilde{N}_++1}}\, .
$$
A first integration between $r$ and $2r$, with $r>0$ small enough, yields
$$
-r^{\tilde{N}_+-1}(u-\gamma v)'(r)> (2r)^{\tilde{N}_+-1}(u-\gamma v)'(2r)-r^{\tilde{N}_+-1}(u-\gamma v)'(r)> c_4 r^{\tilde{N}_+-\frac{2p}{p-1}}\, ,
$$
with $c_4=\frac{c_3}{\Lambda} \frac{2^{\tilde{N}_+-\frac{2p}{p-1}}-1}{\tilde{N}_+-\frac{2p}{p-1}}>0$.  By integrating once more between $r$ and $2r$, we obtain
$$
(u-\gamma v)(r)> (u-\gamma v)(r)- (u-\gamma v)(2r)>\frac{c_5}{r^{\frac{2}{p-1}}}\, ,
$$
with $c_5= c_4 \frac{p-1}{2}\left( 1- 2^{-\frac{2}{p-1}}\right)>0$. On the other hand, by the growth assumption on $v$ and the definition of $\gamma'$ we also have
$$
(u-\gamma v)(r)=v(r) \left( \frac{u(r)}{v(r)}-\gamma\right)\leq \frac{c_2 (\gamma'-\gamma)}{r^{\frac{2}{p-1}}}\, ,
$$
so that we conclude
$$
c_2 (\gamma'-\gamma) > c_5\, .
$$
As $\gamma\to \gamma'$, the left hand side tends to zero, whereas the right hand side becomes
$$
c_5=c_5(\gamma)\to \frac{p-1}{2}\left( 1- 2^{-\frac{2}{p-1}}\right)\frac{2^{\tilde{N}_+-\frac{2p}{p-1}}-1}{\tilde{N}_+-\frac{2p}{p-1}}   \frac{((\gamma')^p-\gamma') c_1^p}{2 \Lambda}>0\, .
$$
This yields a contradiction for $\gamma$ sufficiently close to $\gamma'$. Thus, we have   $\gamma'\leq 1$, that is $u(r)\leq v(r)$.

 \end{proof}

 \section{Proof of  Theorem \ref{mainresult} } 
Let us start with an explicit computation producing sub- and super- solutions: for any $c>0$, $\gamma >0$, we observe that the function $z(r) =c r^{- \gamma}$ satisfies
\begin{equation}\label{gamma}
{\cal M}^+(D^2z)+\mu \frac{z}{r^2}=c\, \Lambda\left (\gamma^2-(\tilde{N}_+-2)\gamma+\frac{\mu}{\Lambda}\right)r^{- \gamma-2}= c\, \Lambda ( \gamma-\tau^+) ( \gamma-\tau^-)r^{- \gamma-2}\, .
\end{equation}

Hence, $z(r) =c r^{- \gamma}$ are super solutions of \eqref{eqmain} for any $\gamma\in [\tau^-,\tau^+]$ and any $c>0$, independently of $p$.
On the other hand, for $\gamma \in  (0, \tau^-)\cup(\tau^+,+\infty)$, one has to relate the values of $\gamma,\ p$ and $c$. In particular, for $\gamma=\frac{2}{p-1}\in  (0, \tau^-)\cup(\tau^+,+\infty)$  the functions
$z=c r^{-2\over p-1} $ are respectively sub-solutions and  super-solutions if 
\begin{equation}\label{subsuper}
c \leq K \, ,  \mbox{ respectively } c \geq K ,
\end{equation}
with $K: = \left[\Lambda \left( \frac{2}{p-1}-\tau^+\right) \left( \frac{2}{p-1}-\tau^-\right)\right]^{\frac{1}{p-1}}$.

Hence, for either  $p<p^\star$ or $p>p^{\star \star}$,  
$u(r)=Kr^{-2\over p-1} $ is an explicit solution. 
This  proves statement \eqref{2surp-1} of case 1 and case 4.
\medskip

Let us now concentrate on the remaining cases for $p<p^\star$. 

In order to show the existence of a solution $v$
 asymptotic to $r^{-\tau^+}$ near zero, we 
 will use an approximation argument combined with  the   comparison theorem given by Proposition \ref{compa1}.                                
      
 For $\epsilon>0$ and $0<\delta<\tau^+$, let  
 $$\underline{v} (r)= \epsilon (r^{-\tau^+}- r^{-\tau^++\delta})\, .$$
We observe that in $B_1\setminus \{0\}$ one has $\underline{v}>0\, ,\ \underline{v}^\prime < 0$,  $\underline{v}^{\prime \prime}>0$ and
 \begin{eqnarray*}
                                   {\cal M}^+( D^2  \underline{v} )+ \mu {\underline{ v} \over r^2} &=& \epsilon r^{-\tau^++\delta-2} \Lambda \delta ( 2(\tau^+-\tau)-\delta)\, .
 \end{eqnarray*}
 By choosing $0<\delta< \min\{ \tau^+\, ,\ 2-(p-1)\tau^+\, ,\ 2(\tau^+-\tau)\}$ and $0<\epsilon^{p-1} \leq \Lambda \delta ( 2(\tau^+-\tau)-\delta)$,
 we  easily  obtain that $\underline{v} $ is a sub-solution of \eqref{eqmain}. 
 
 More in general, we observe that for $0<\delta< \min\{ \tau^+\, ,\ 2-(p-1)\tau^+\, ,\ 2(\tau^+-\tau)\}$ and for any $a>0$, the function
 $$
 \underline{z} (r) = a\, r^{-\tau^+}- \frac{a^p}{\Lambda \delta \left( 2(\tau^+-\tau)-\delta\right)} r^{-\tau^++\delta}
 $$
 is a subsolution in $B_{r_0}\setminus \{0\}$ for $r_0^\delta\leq \frac{\Lambda \delta \left( 2(\tau^+-\tau)-\delta\right)}{a^{p-1}}$.

 On the other hand, by the initial computation, we also have that $\overline{v}(r)=r^{-\tau^+}$ is a super-solution of \eqref{eqmain}.
 
 Let us now consider any decreasing sequence $\{ r_n\} \subset (0,1)$, and let us  define a sequence of approximating radial solutions $\{v_n\}$, with $v_n: [r_n,1]\to [0, +\infty)$, in the following way: $v_1:[r_1,1]\to [0,+\infty)$ is the solution of the Dirichlet problem
 $$
 \left\{
 \begin{array}{c}
\displaystyle -\mathcal{M}^+(D^2v_1) + |v_1|^{p-1}v_1= \mu  \frac{\underline{v}}{r^2} \quad \hbox{ in } B_1\setminus \overline{B}_{r_1}\\[1ex]
 v_1(1)=0\, ,\ v_1(r_1)= \underline{v}(r_1)
 \end{array}
 \right.
 $$
 and, inductively, $v_{n+1}:[r_{n+1},1]\to [0,+\infty)$ is the solution of
 $$
 \left\{
 \begin{array}{c}
\displaystyle -\mathcal{M}^+(D^2v_{n+1}) + |v_{n+1}|^{p-1}v_{n+1} =   \frac{\mu}{r^2}\left\{
\begin{array}{ll}
v_n & \hbox{ in } B_1\setminus B_{r_n}\\[1ex]
\underline{v} & \hbox{ in } B_{r_n}\setminus \overline{B}_{r_{n+1}}
\end{array} \right. \\[3ex]
 v_{n+1}(1)=0\, ,\ v_{n+1}(r_{n+1})= \underline{v}(r_{n+1})
 \end{array}
 \right.
 $$
  By using induction and the standard comparison principle for the operator $-\mathcal{M}^+(\cdot) +|\cdot|^{p-1} (\cdot)$, it easily follows that $v_n$ is radial, positive for $r<1$ and it satisfies $v_n\geq \underline{v}$ in $[r_n,1]$ for every $n\geq 1$. This in turn implies that $v_{n+1}\geq v_n$ in $[r_n,1]$ for all $n\geq 1$. 
 Moreover, Proposition \ref{compa1} yields that $v_n\leq \overline{v}$ in $[r_n,1]$ for all $n\geq 1$. Standard stability and regularity results then imply that the monotone locally uniformly bounded sequence $\{v_n\}$ converges in $C^2_{{\rm loc}}((0,1])$ to a radial solution $v$ of equation \eqref{eqmain}  satisfying 
 $$
 \underline{v}\leq v \leq \overline{v} \quad \hbox{ for } 0<r\leq 1\, .$$
It is  left to prove that $v(r)$ is asymptotic to $r^{-\tau^+}$ as $r\to 0$. 
By the bounds above,  there exists   a decreasing sequence $\{ r_n\}$ converging  to zero  such that  $v(r_n) r_n^{\tau^+} \rightarrow c_1$ for some $c_1>0$. 
For $\epsilon >0$ arbitrarily small, let $\nu \geq 1$ such that
  $ c_1-\epsilon \leq  v(r_n) r_n^{\tau^+} \leq c_1+\epsilon$ for $n\geq \nu$ and, moreover,  $r_\nu^\delta <\frac{\Lambda \delta \left( 2(\tau^+-\tau)-\delta\right)}{(c_1-\epsilon)^{p-1}}$.  
  
Then, for all $n\geq \nu$, we can consider on the interval $[r_{n+1}, r_n]$ the supersolution
     $\overline{z}  = (c_1+\epsilon) r^{-\tau^+}$ and the subsolution 
     $\underline{z}= (c_1-\epsilon)  r^{-\tau^+}- {(c_1-\epsilon)^{p} \over \Lambda \delta(2 ( \tau^+-\tau)-\delta )} r^{-\tau^++\delta }$. 
 Since     
 $$ \underline{z}(r_n) \leq (c_1-\epsilon ) r_n^{-\tau^+} \leq v(r_n)\leq (c_1+\epsilon) r_n^{-\tau^+}= \overline{z}(r_n)$$
 as well as 
     $$ \underline{z}(r_{n+1} )\leq (c_1-\epsilon ) r_{n+1}^{-\tau^+} \leq v(r_{n+1})\leq (c_1+\epsilon) r_{n+1}^{-\tau^+}= \overline{z}(r_{n+1})\, ,$$
       the comparison principle given by Proposition \ref{compa1} yields 
       $\underline{z}\leq v \leq  \overline{z} $ in $[r_{n+1}, r_n]$ for all $n\geq \nu$, which implies \eqref{tau+}.
\medskip                                          
                                          
In the case $p<p^{\star \star}=1+\frac{2}{\tau^-}$, an analogous construction can be repeated by replacing $\tau^+$ with $\tau^-$. In particular, the function
$$  \underline{w}(r)= \epsilon \left( r^{-\tau^-} + r^{-\tau^-+\delta}\right)$$
 is a sub-solution in $B_1\setminus \{0\}$ provided that  $0< \delta < \min \{ \tau^-\, ,\ 2-(p-1)\tau^-\}$ and $\epsilon^{p-1}\leq \frac{\Lambda \delta \left( 2(\tau-\tau^-)+\delta\right)}{2^p}$. More in general, 
 $ \underline{w}$   is a sub-solution   in    $B_{r_0}\setminus \{0\}$   for any   $0< \delta < \min \{ \tau^-\, ,\ 2-(p-1)\tau^-\}$, $\epsilon >0$ and 
 $$r_0^{2-(p-1)\tau^--\delta}\leq \frac{\Lambda \delta \left( 2(\tau-\tau^-)+\delta\right)}{2^p \epsilon^{p-1}}\, .$$             
By using also the supersolution $\overline{w} (r)=c r^{-\tau^-}$ for suitable  $c > 0$ and repeating the approximating procedure we performed above, we obtain the existence of a radial solution $w$ satisfying
$$
\underline{w}(r) \leq w(r)\leq \overline{w}(r)\quad \hbox{ for } 0<r\leq 1\, .
$$
The solution $w$ is further showed to satisfy \eqref{tau-}, thus yielding the proof of   statement \eqref{tau-} in  1  and  of statement  2.                                                                         
\medskip                                          
                                         
Let us finally prove statement 3. For  $p=p^{\star \star}$, i.e. $\frac{2}{p-1}=\tau^-$, let us                                          
 consider the function
                                            $$ \underline{u} (r) = a\,  r^{-\tau^-} ( -\log r+c)^{-\tau^-\over 2}\, ,$$
 with $a>0$ and $c>1$.  Then, a direct computation shows that $\underline{u}'<0 ,\ \underline{u}''>0$ and
 $$
\mathcal{M}^+(D^2\underline{u}) + \mu \frac{\underline{u}}{r^2}=
\frac{\Lambda a  \tau^-}{r^{\tau^-+2}}\left( \frac{\tau-\tau^-}{(-\log r+c)^{\tau^-/2 +1}}+ \frac{\tau^-+2}{4(-\log r+c)^{\tau^-/2 +2}}\right)
$$
Since
$$
\underline{u}^p(r)= \frac{a^p}{r^{\tau^-+2}(-\log r+c)^{\tau^-/2 +1}}\, ,
$$
we easily deduce that $\underline{u}$ is a sub-solution for  $a^{p-1}\leq \Lambda \tau^- (\tau -\tau^-)=\bar K^{p-1}$.

On the other hand,   an analogous computation shows that, for any $b\, ,\ \delta>0$  and $c>1$ the function
$$\overline{u} (r)= b\left( r^{-\tau^-}  ( -\log r+c)^{-\tau^-/2} -r^{-\tau^-} ( -\log r+c)^{-\tau^-/2 -\delta}\right)$$
 satisfies
 $$
 \mathcal{M}^+(D^2\overline{u}) + \mu \frac{\overline{u}}{r^2}\leq 0\,.
 $$
 Hence, by   using the subsolution
 $\underline{u}$ with $a=\bar K$
 and the supersolution $\overline{u}$ for any fixed $b>\bar K$ and $c^\delta \geq \frac{b}{b-\bar K}$, we can repeat the approximation procedure we performed above in order to   obtain a solution $u$ satisfying
 $$
 \underline{u}\leq u\leq \overline{u} \quad \hbox{ in } B_1\setminus\{0\}\, .
 $$
 This implies
 $$
 \bar K\leq \liminf_{r\to 0} r^{\tau^-} (-\log r)^{\tau^-/2}u(r) \leq \limsup_{r\to 0} r^{\tau^-} (-\log r)^{\tau^-/2}u(r)\leq b
 $$
 for any $b>\bar K$, that is statement 3.
 
\hfill $\Box$

\section{Proof of Theorem \ref{behave}}
This section will be devoted to the proof of Theorem \ref{behave}, after some preliminary observations.

First of all, we note that if $u$ is a nonnegative solution of equation \eqref{eqmain} then, by the strong minimum principle, either $u$ vanishes identically or it is strictly positive in $B_1\setminus \{0\}$. Thus, in the sequel, by a solution $u$ we mean a strictly positive radial solution.

\begin{prop}\label{bounded}
Let  $u$ be a solution of \eqref{eqmain}. Then
\begin{equation}\label{uunbounded}
\limsup_{r\to 0} u(r)=+\infty
\end{equation}
and
\begin{equation}\label{xbounded}
\limsup_{r\to 0} u(r) r^{2\over p-1}< +\infty\,.
\end{equation}
\end{prop}        
        
 \begin{proof}
In order to prove \eqref{uunbounded}, let us argue by contradiction and assume that  $u$ is bounded. Then, for $r$ sufficiently small, one has
         $${ \cal M}^+ ( D^2 u )= u^p- \mu \frac{u}{r^2}\leq  0\, .$$ 
This implies in particular that,    whatever  the signs of $u^\prime$ and $u^{ \prime \prime}$ are,  
          $$ \Lambda u^{ \prime \prime} + \lambda ( N-1) \frac{u^\prime}{ r} \leq \mathcal{M}^+(D^2u)\leq 0\, ,$$
that is
          $$ (u^\prime r^{ \tilde{N}_+-1})^\prime \leq 0\, .$$
Hence, the limit $ c\, :=\lim_{r\to 0}u^\prime(r) r^{\tilde{N}_+-1}$ exists. If $c>0$, one gets that  $u$ is negative large, whereas if $c<0$, it follows that $u$ is unbounded from above. Thus, we deduce that $ c=0$, so that 
 $u^\prime (r)\leq 0$ and   $0<u(r)\leq u(0)$ for any $0<r<1$.
 By further observing that  $u^p-\mu \frac{u}{ r^2} \leq -\frac{ \mu}{2} \frac{u}{ r^2}$ for $r$ small enough, we get
         $$ (u^\prime (r) r^{\tilde{N}_+-1})^\prime \leq -  u(0) \frac{ \mu}{4}  r^{ \tilde{N}_+-3}\, .$$
Hence, for some $c_2>0$,
         $$ u^\prime (r)\leq - \frac{c_2}{r}\, ,$$
          which yields
          $$ u(r) \geq -c_2\log r + c_3\, ,$$ a contradiction to the boundedness of $u$.  

\medskip
Let us now prove \eqref{xbounded}. We preliminarily observe that if $u(r)r^{\frac{2}{p-1}}$ is bounded along a decreasing sequence $\{r_n\}$ converging to zero, then it is bounded in a neighborhood of zero. Indeed, let us assume that there exists $C>0$ such that 
$$
u(r_n)r_n^{\frac{2}{p-1}}\leq C\, .
$$
By using the computation made at the beginning of the previous section, we know that, for $C>0$ sufficiently large, the function 
 $C r^{-2\over p-1}$ is a super-solution of  equation \eqref{eqmain}. Then, by applying Theorem 
 \ref{compa1} in every  annulus  $B_{r_{n}}\setminus\overline{B}_{r_{n+1}}$, we obtain  
$$u(r)\leq C r^{-2\over p-1} \quad \hbox{ for } 0<r\leq r_1\, .
$$ 
Hence, arguing  by contradiction  again and assuming that \eqref{xbounded} fails, we deduce
$$
\lim_{r\to 0} u(r) r^{\frac{2}{p-1}}= +\infty\, .
$$    
This implies that, for $r$ small enough, one has
$$
\mathcal{M}^+(D^2u)= u^p -\mu \frac{u}{r^2} \geq { u^p \over 2}>0\, .$$ 
It then follows that $u'$ has constant sign in a right neighborhood of zero. Indeed, if not, we could find $0<s<t<1$ such that $u'(s)=u'(t)=0$ and $u'(r)<0$ for $s<r<t$. Then, there would exist some $\sigma \in (s,t)$ such that $u''(\sigma)=0$, so that
$$
\mathcal{M}^+(D^2u (\sigma))<0\, .$$
By using also \eqref{uunbounded}, we deduce $
u^\prime \leq 0$. We then have 
$$ u^{ \prime \prime} \geq { u^p\over 2 \Lambda}\, ,$$
 and then
          $$u^{ \prime \prime} u^\prime \leq {u^p u^\prime \over 2 \Lambda}\, .$$
           Integrating, one gets 
           $$ (u^\prime )^2(r) - c  u^{p+1} (r) \geq  (u^\prime )^2(r_0) - c  u^{p+1} (r_0)$$
            for $r < r_0$ and for some $c>0$. Using $u^{p+1}(r)\rightarrow +\infty$ as $r\to 0$, for $r$ small enough we get
            $$(u^\prime)^2 \geq \frac{c}{2} u^{p+1}$$
            or, equivalently,
            $$ {-u^\prime \over u^{ p+1\over 2}} \geq c_1>0\, .$$
This yields
             $$ u^{ 1-p \over 2} \geq c_2 r\, ,$$
 a contradiction to the assumption. 
  
\end{proof}

Other general asymptotic properties of any solution $u$ are given in the next results.

\begin{prop}\label{signeuprime} Let $u$ be a solution of \eqref{eqmain}. Then:
 \begin{enumerate}
 
 \item If $\dis\lim_{r\rightarrow 0}u(r) r^{2\over p-1}=0$, then 
  $u^\prime (r)\leq 0$ and   $u^{ \prime \prime}(r) \geq 0$ for $r$ sufficiently small;

\item if $\dis\limsup_{r\rightarrow 0} u(r) r^{2\over p-1}>0$ and either $1<p<p^\star$ or $p>p^{\star \star}$,  then
        $$\lim_{r\to 0} u(r) r^{\frac{2}{p-1}}= \left[ \Lambda \left( {2\over p-1}-\tau^+\right) \left( {2\over p-1}-\tau^-\right)\right]^{1\over p-1}= K\, .$$

        \end{enumerate}
\end{prop}

\begin{proof} 
Let $u$ be a solution satisfying $\lim_{r\to 0} u(r)r^{\frac{2}{p-1}}=0$. Then, for $r$ small enough, one has
$$
\mathcal{M}^+(D^2u)= \frac{u}{r^2} \left( u^{p-1}r^2-\mu\right)<0 \, .
$$
By arguing as in the proof of Proposition \ref{bounded}, we deduce that $u'$ has constant sign for $r>0$ sufficiently small, whence, by \eqref{uunbounded}, we get $u'\leq 0$.

Furthermore, we have that $u''(r)$ also has constant sign for $r>0$ small enough. Indeed, if not, we can find a decreasing sequence $\{ r_n\}$ converging to zero such that  for all $n\geq 1$ one has $u''(r_n)=0$ and 
$$
\begin{array}{rl}
u''(r) (r-r_{2n+1})> 0 & \hbox{ for $r$ close to $r_{2n+1}$, $r\neq r_{2n+1}$}\, ,\\[1ex]
u''(r)(r-r_{2n})< 0  &  \hbox{ for $r$ close to $r_{2n}$, $r\neq r_{2n}$}\, .
\end{array}
$$
On the other hand, for $r>0$ sufficiently small, we can consider the function
$$\varphi ( r) = -(\tilde{N}_+-1) r u^\prime (r)- \frac{\mu}{\Lambda}  u (r) +  \frac{r^2}{\Lambda} u^p (r)\, .$$
We  observe that $\varphi(r)$ coincides with $u^{\prime\prime}(r)$ up to a positive factor, so that they  share the same sign. 
Moreover, one has
$$ \varphi^\prime (r) = -(\tilde{N}_+-1) u'(r) -(\tilde{N}_+-1) r u''(r) - \frac{\mu}{\Lambda} u^\prime(r) +\frac{2 r}{\Lambda} u^p(r)+ p \frac{r^2 }{\Lambda} u^{p-1}(r) u^\prime (r)\, .
$$
By evaluating the above derivative for $r=r_n$ and using that $u''(r_n)=0$, $u'(r_n)<0$ and 
$$\lim_{n\to +\infty} u^{p-1}(r_n)r_n^2=0\, ,$$ 
we obtain
 $$ \varphi^\prime ( r_n)\geq \frac{2r_n}{\Lambda} u^p(r_n) >0\, .$$
This implies that  $\varphi(r)(r-r_n) >0$ for $r$ close to $r_n$, $r\neq r_n$, for $n$ large.  Hence $u''(r)(r-r_n) >0$ for $r$ close to $r_n$, $r\neq r_n$, for all $n$ sufficiently large, contradicting the choice of $r_n$. 

Thus, we obtain that $u''(r)$ does not change sign for $r$ small enough. Since a positive, concave and decreasing function is necessarily bounded on a bounded interval, by \eqref{uunbounded} we deduce $u''(r)\geq 0$ for $r>0$ small enough.        
\medskip

In order to prove statement 2., let us set $l\, := \limsup_{r\rightarrow 0} u(r) r^{2\over p-1}$, which is supposed to be positive. Moreover, by   Proposition \ref{bounded}, we have  that $l<+\infty$. Let us   select a decreasing sequence $\{ r_n\}$ converging to zero such that
$$
u(r_n)r_n^{\frac{2}{p-1}}\geq \frac{l}{2} \quad \hbox{ for all } n\geq 1\, .
$$
By the assumption on $p$ and by \eqref{subsuper}, we know that the function $z(r)=\min\left\{K, \frac{l}{2}\right\} r^{-\frac{2}{p-1}}$ is a sub-solution of equation \eqref{eqmain}. Applying Theorem \ref{compa1} in each annulus $B_{r_n}\setminus \overline{B}_{r_{n+1}}$, we then deduce
$$
u(r)r^{\frac{2}{p-1}} \geq \min\left\{K, \frac{l}{2}\right\} >0\quad \hbox{ for } 0<r\leq r_1\, .
$$
In order to get the conclusion of statement 2., we are going to use   Theorem \ref{compa2} applied with  suitable constructed sub- and super-solutions. We will prove  both that $\liminf_{r\to 0} u(r)r^{\frac{2}{p-1}}\geq K$ and that $\limsup_{r\to 0} u(r)r^{\frac{2}{p-1}}\leq K$.

Arguing by contradiction, let us assume that $l_1\, :=\liminf_{r\to 0} u(r)r^{\frac{2}{p-1}}< K$ and let us select $K_1>0$ such that
$$
l_1<K_1<K\, .
$$
Let us further fix $r_0>0$ suitably small such that $u(r_0)  r_0^{\frac{2}{p-1}}<K_1$ and,
for $a>0$ and $0<\gamma <\frac{2}{p-1}$ to be  fixed, let us consider the function
$$
\underline{u}(r)= \frac{K_1}{r^{\frac{2}{p-1}}} -\frac{a}{r^\gamma}\quad \hbox{ for } 0<r\leq r_0\, .
$$
We observe that, by choosing the parameters  $a$ and $\gamma$ such that
$$
K_1-u(r_0) r_0^{\frac{2}{p-1}}\leq a\, r_0^{\frac{2}{p-1}-\gamma}< K_1\, ,
$$
 
then $\underline{u}(r)$ is a convex decreasing function for $0<r\leq r_0$, satisfying also
$$
\underline{u}(r_0)\leq u(r_0)\, .
$$
Moreover, an analogous computation as for \eqref{gamma} yields
$$
\mathcal{M}^+(D^2\underline{u})+\mu \frac{\underline{u}}{r^2}=
\frac{K_1 K^{p-1} - a \Lambda (\gamma-\tau^-)(\gamma-\tau^+) r^{\frac{2}{p-1}-\gamma}}{r^{\frac{2p}{p-1}}}\, .
$$
By selecting $\gamma$ sufficiently close to $\frac{2}{p-1}$ and by
further imposing that
$$
a\, r_0^{\frac{2}{p-1}-\gamma}\leq \frac{K_1 \left( K^{p-1}-K_1^{p-1}\right)}{\Lambda (\gamma-\tau^-)(\gamma-\tau^+)}\, ,
$$
it then follows that
$$
\mathcal{M}^+(D^2\underline{u})+\mu \frac{\underline{u}}{r^2}\geq \frac{K_1^p}{r^{\frac{2p}{p-1}}}\geq \underline{u}^p\quad \hbox{ in } B_{r_0}\setminus\{0\}\, .
$$
We note that all the imposed requirements are satisfied by setting
$$
a=\frac{K_1-u(r_0) r_0^{\frac{2}{p-1}}}{r_0^{\frac{2}{p-1}-\gamma}}
$$
and by choosing $K_1$ and $r_0$   such that $K_1-u(r_0) r_0^{\frac{2}{p-1}}>0$ is sufficiently small.

Hence, we have obtained that $\underline{u}$ is a subsolution of equation \eqref{eqmain} in $B_{r_0}\setminus\{0\}$, satisfying $\underline{u}(r_0)=u(r_0)$. Theorem \ref{compa2} then yields 
$$
u(r)\geq \underline{u}(r)\quad \hbox{ for }0<r\leq r_0\, ,
$$
which implies the contradiction
$$
\liminf_{r\to 0} u(r)r^{\frac{2}{p-1}} =l_1\geq K_1\, .
$$
The reached contradiction proves that
$$
\liminf_{r\to 0} u(r)r^{\frac{2}{p-1}} =l_1\geq K\, .
$$
Analogously, by assuming that $l=\limsup_{r\to 0}  u(r)r^{\frac{2}{p-1}}>K$ and by fixing $K_2>0$ such that
$$
K<K_2< l
$$
and $r_0>0$ suitably small that $K_2< u(r_0)r_0^{\frac{2}{p-1}}$, one can find $a>0$ and $0<\gamma <\frac{2}{p-1}$ such that the function
$$
\overline{u}(r)= \frac{K_2}{r^{\frac{2}{p-1}}} +\frac{a}{r^\gamma}\quad \hbox{ for } 0<r\leq r_0
$$
is a super-solution of \eqref{eqmain} in $B_{r_0}\setminus \{0\}$ satisfying $u(r_0)\leq \overline{u}(r_0)$. Theorem \ref{compa2} then gives
$$
u(r)\leq \overline{u}(r)\quad \hbox{ for } 0<r\leq r_0\, ,
$$
and, therefore,
$$
l=\limsup_{r\to 0} u(r)r^{\frac{2}{p-1}} \leq K_2\, .
$$
The reached contradiction shows that
$$
l=\limsup_{r\to 0} u(r)r^{\frac{2}{p-1}} \leq K\, ,
$$
which completely proves statement 2.

 \end{proof} 
 
\begin{prop}\label{monotonicity}
Let $p^\star \leq p\leq p^{\star \star}$. Then, the function
$$
v(r)=r^{\tau^-} u(r)
$$
is monotone non decreasing for $0<r<1$. As a consequence,  $u(1)>0$.
\end{prop} 

\begin{proof}
 Let $0<r_0<1$ be arbitrarily fixed and, for $0<\epsilon< u(r_0)r_0^{\frac{2}{p-1}}$,  let us consider the function
$$
\overline{u}_\epsilon(r)= \frac{\epsilon}{ r^{2\over p-1}} + \frac{c_\epsilon}{r^{\tau^-}}\, ,$$
with $c_\epsilon=\frac{u(r_0)r_0^{\frac{2}{p-1}}-\epsilon}{r_0^{\frac{2}{p-1}-\tau^-}}$.
  Then, one has $\overline{u}_\epsilon(r_0)=u(r_0)$ and, by direct computation,
  $$
  \mathcal{M}^+ ( D^2 \overline{u}_\epsilon) + \frac{\mu}{r^2} \overline{u}_\epsilon= \epsilon \Lambda \frac{\left( \frac{2}{p-1}-\tau^-\right)\left( \frac{2}{p-1}-\tau^+\right)}{r^{\frac{2p}{p-1}}}\leq 0\, .$$
Hence, $\overline{u}_\epsilon$ is trivially a super-solution of \eqref{eqmain}, so that,  by the comparison Theorem \ref{compa2}, one gets 
$$ u(r) \leq \overline{u}_\epsilon(r)\quad \hbox{ for all } 0<r\leq r_0\, .$$
By letting $\epsilon\to 0$, we deduce
$$
u(r)\leq \frac{u(r_0)r_0^{\tau^-}}{r^{\tau^-}} \quad \hbox{ for all } 0<r\leq r_0\, ,
$$
that is the conclusion.

 \end{proof}                  
 
 Our last preliminary result concerns the case $p=p^{\star \star}$. It provides a first rough upper bound which will be used in the proof of the sharp result given by statement 3 of Theorem \ref{behave}.
 
 \begin{prop} \label{p**}
If $p = p^{\star \star}$, then there exists a  positive constant $C$  such that
$$ u(r) \leq  \frac{C}{ \left[(- \log r )^{1\over 2} r\right]^{\tau^-}}\, .$$
In particular, one has  $\lim_{r \to 0} u (r) r^{2\over p-1} = 0$. 
 \end{prop}

 \begin{proof}
 We recall that $p=p^{\star \star}$ means that $\frac{2}{p-1}=\tau^-$.  
 
For $\epsilon>0$ arbitrarily small and  $C>0$ to be suitably fixed, let us  introduce the function 
$$  \overline{u}(r)= \frac{C}{ \left[ (-\log ( \epsilon + r))^{1\over 2} r\right]^{\tau^-}}\, .$$
We observe that there exists $r_0>0$, independent of $\epsilon$, such that $\overline{u}'(r)\leq 0$ and $\overline{u}''(r)\geq 0$ for $0<r\leq r_0$.  In particular, we fix $r_0$ such that 
$$
-\log (r_0+\epsilon)\geq \frac{1}{2}\, .
$$
Then, by directly  computing and dropping the negative terms, we obtain
$$
\begin{array}{ll}
\displaystyle \mathcal{M}^+(D^2\overline{u})+\mu \frac{\overline{u}}{r^2}  &  \displaystyle \leq \Lambda C \left[ \frac{\tau^- (\tilde{N}_+-1)}{2r^{\tau^-+1}(r+\epsilon) (-\log (r+\epsilon))^{\frac{\tau^-}{2}+1}}+
\frac{\tau^- (\tau^-+2)}{4 r^{\tau^-}(r+\epsilon)^2(-\log (r+\epsilon))^{\frac{\tau^-}{2}+2}}\right] \\[4ex]
& \displaystyle \leq \frac{\Lambda}{2} C \frac{\tau^- ( \tau^- +\tilde{N}_++1)}{r^{\tau^-+2}(-\log (r+\epsilon))^{\frac{\tau^-}{2}+1}}
\end{array}
$$
By the assumption on $p$, one has $\tau^-+2=p \tau^-$ and $\frac{\tau^-}{2}+1 = p\frac{\tau^-}{2}$. Hence, we deduce
$$
\mathcal{M}^+(D^2\overline{u})+\mu \frac{\overline{u}}{r^2}\leq \overline{u}^p= \frac{C^p}{\left[ (-\log ( \epsilon + r))^{1\over 2} r\right]^{p\tau^-}}
$$
provided that
$$C^{p-1} >\frac{\Lambda}{2}  \tau^- ( \tau^- +\tilde{N}_++1)\, .$$ 
Moreover, we notice that, for $0<r\leq r_0$, one has
$$
\frac{C(-\log \epsilon)^{\frac{-\tau^-}{2}}}{ r^{\frac{2}{p-1}}}\leq \overline{u}(r)\leq \frac{2^{\frac{\tau^-}{2}} C}{r^{\frac{2}{p-1}}}\, .
$$
Hence, by further imposing that
$$
u(r_0)\leq \overline{u}(r_0)=\frac{C}{ \left[ (-\log ( \epsilon + r_0))^{1\over 2} r_0\right]^{\tau^-}}\, ,
$$
that is
$$
C\geq u(r_0) \left[ (-\log ( \epsilon + r_0))^{1\over 2} r_0\right]^{\tau^-}\, ,
$$
 we can apply Theorem \ref{compa2}, which yields
 $$
 u(r)\leq \overline{u} (r)\quad \hbox{ for } 0<r\leq r_0\, .
 $$
 By letting $\epsilon \to 0$, we get the conclusion.
 
 \end{proof}

At this point, we have all the tools to be used in the proof of Theorem \ref{behave}.

 \subsection{Proof of 1. }
 
 Let $u$ be a solution of \eqref{eqmain}, with $1<p<p^\star$. 
 
 Two cases are possible: either $\limsup_{r\to 0} u(r)r^{\frac{2}{p-1}}>0$ or $\lim_{r\to 0} u(r)r^{\frac{2}{p-1}}=0$.  
 
 In the first case,  by statement 2 of Proposition \ref{signeuprime}, we get that $u$ satisfies \eqref{2surp-1}. 
 
 In the second case,  by Proposition \ref{signeuprime} again, we obtain that $u(r)$ is convex and decreasing for $r$ sufficiently small. Thus, we can write more explicitly the equation satisfied by $u$, that is
 \begin{equation}\label{equ}
 u'' +(\tilde{N}_+-1) \frac{u'}{r} +\frac{\mu}{\Lambda} \frac{u}{r^2} - \frac{u^p}{\Lambda}=0\, .
 \end{equation}
Let us introduce the so called  Emden--Fowler transformed of $u$, namely the function
 $$ x(t) = e^{\frac{2}{p-1} t} u(e^t)\qquad \hbox{ for } t\leq 0\, .$$
 A direct computation shows that equation \eqref{equ} translates in terms of $x(t)$ into the autonomous equation
 \begin{equation}\label{eqx}
 x''- 2 \left( \frac{2}{p-1} -\tau\right) x' +\left(  \frac{2}{p-1} -\tau^+\right) \left( \frac{2}{p-1}-\tau^-\right)x -\frac{x^p}{\Lambda}=0\, .
\end{equation}
Thus, we have that $x(t)$ is a solution of \eqref{eqx}, defined for $t\in (-\infty, T]$ for some $T\leq 0$, and satisfying 
$$
\lim_{t\to -\infty} x(t)=0\, .
$$ 
Let us observe that equation \eqref{eqx} can be written in the form
\begin{equation}\label{eqxlambda}
x''-\left( \lambda_1+\lambda_2\right) x'+\lambda_1 \lambda_2 \, x= \frac{x^p}{\Lambda}\, ,
\end{equation}
where
\begin{equation}\label{l1l2}
\begin{array}{c}
\lambda_1\, :=\frac{2}{p-1}-\tau^+\\[2ex]
\lambda_2\, :=\frac{2}{p-1}-\tau^-
\end{array}
\end{equation}
are the roots of the characteristic equation associated with equation \eqref{eqx} linearized at zero. 

We note that, by the assumption on $p$, one has
$$
0<\lambda_1<\lambda_2\, .
$$
By the  theory for perturbed linear dynamical systems (see e.g. Theorem 3.5, case (ii), of \cite{Har}) applied to the  positive solution $x(t)$  converging to zero as $t\to -\infty$, one has that either there exists $c_1>0$ such that
\begin{equation}\label{asym1}
x(t) e^{-\lambda_1 t} \to c_1 \quad \hbox{ as } t\to -\infty\, ,
\end{equation}
or 
there exists $c_2>0$ such that
\begin{equation}\label{asym2}
x(t) e^{-\lambda_2 t} \to c_2 \quad \hbox{ as } t\to -\infty\, .
\end{equation}
By the definition of $x(t)$ and of $\lambda_1$ and $\lambda_2$, this amounts to say that $u$ satisfies either \eqref{tau+} or \eqref{tau-}, i.e. the remaining cases stated in 1 of  Theorem \ref{behave}.

For the sake of completeness, let us provide a direct proof of \eqref{asym1} and \eqref{asym2}.

Since $x(t)$ is a solution of \eqref{eqxlambda} in $(-\infty, T]$, then $x(t)$ satisfies the integral formulation
\begin{equation}\label{integral}
x(t)=x_1 e^{\lambda_1 (t-T)}+x_2 e^{\lambda_2 (t-T)} +\frac{e^{\lambda_1 t}}{\Lambda\, (\lambda_2-\lambda_1)}\int_t^T x(s)^p e^{-\lambda_1 s}\, ds -\frac{e^{\lambda_2 t}}{\Lambda\, (\lambda_2-\lambda_1)}\int_t^T x(s)^p e^{-\lambda_2 s}\, ds\, ,
\end{equation}
with $x_1\, ,\ x_2\in \R$ related to the initial conditions by the system
$$
\left\{
\begin{array}{c}
x_1+x_2=x(T)\\[1ex]
\lambda_1\, x_1+\lambda_2\, x_2= x'(T)
\end{array}
\right.
$$
Moreover, since $x(t)$ is positive and $x(t)\to 0$ as $t\to -\infty$, we can fix $\delta>0$ as small as needed and select $T$ sufficiently negative such that
$$
0<x(t)\leq \delta \qquad \hbox{ for all } t\leq T\, .
$$
From \eqref{integral} it then follows that
$$
x(t) e^{-\lambda_1 t} \leq x_1 e^{-\lambda_1 T}+x_2 e^{-\lambda_2 T} +\frac{\delta^{p-1}}{\Lambda\, (\lambda_2-\lambda_1)}\int_t^T x(s) e^{-\lambda_1 s}\, ds\, .
$$
Gronwall's inequality then yields
$$
x(t) e^{-\lambda_1 t} \leq  C\, e^{-\eta t}\, ,
$$
with $\eta= \frac{\delta^{p-1}}{\Lambda\, (\lambda_2-\lambda_1)}$ and for some $C>0$ depending on $\delta$. This in turn implies that
$$
x(t)^pe^{-\lambda_1 t}\leq C^p e^{( (p-1)\lambda_1-p \eta) t}\, ,
$$
so that $x(t)^pe^{-\lambda_1 t}$ is integrable on $(-\infty, T]$ for $\delta$ so small that $\eta< \frac{p-1}{p} \lambda_1$. Coming back to \eqref{integral},  first we  deduce that
$$
0<x(t) \leq C_1\, e^{\lambda_1 t}
$$
for some  $C_1>0$, and  then we conclude that
$$
x(t) e^{-\lambda_1 t} \to x_1 e^{-\lambda_1 T}+ \frac{1}{\Lambda\, (\lambda_2-\lambda_1)}\int_{-\infty}^T x(s)^p e^{-\lambda_1 s}\, ds \, =: c_1\quad \hbox{ as } t\to -\infty\, .
$$
If $c_1>0$, then we get \eqref{asym1}. If, on the contrary, $c_1=0$, then \eqref{integral} becomes
\begin{equation}\label{integral2}
x(t)=x_2 e^{\lambda_2 (t-T)} -\frac{e^{\lambda_1 t}}{\Lambda\, (\lambda_2-\lambda_1)}\int_{-\infty}^t x(s)^p e^{-\lambda_1 s}\, ds -\frac{e^{\lambda_2 t}}{\Lambda\, (\lambda_2-\lambda_1)}\int_t^T x(s)^p e^{-\lambda_2 s}\, ds\, ,
\end{equation}
which immediately yields
$$
0<x(t) \leq C_2\, e^{\lambda_2 t}\, .
$$
Hence, from \eqref{integral2} it follows that
$$
x(t) e^{-\lambda_2 t} \to x_2 e^{-\lambda_2 T}- \frac{1}{\Lambda\, (\lambda_2-\lambda_1)}\int_{-\infty}^T x(s)^p e^{-\lambda_2 s}\, ds \, =: c_2\quad \hbox{ as } t\to -\infty\, .
$$
We claim that, in this case,  $c_2>0$ so that \eqref{asym2} holds true. Indeed, arguing by contradiction, if $c_2=0$ as well, then \eqref{integral2} becomes
$$
x(t)= -\frac{e^{\lambda_1 t}}{\Lambda\, (\lambda_2-\lambda_1)}\int_{-\infty}^t x(s)^p e^{-\lambda_1 s}\, ds +\frac{e^{\lambda_2 t}}{\Lambda\, (\lambda_2-\lambda_1)}\int_{-\infty}^t x(s)^p e^{-\lambda_2 s}\, ds\, .
$$
Using that $x(t) e^{-\lambda_2 t} \to 0$, then we  have, for $t$ sufficiently negative,
$$
x(t) e^{-\lambda_2 t} \leq \frac{1}{\Lambda\, (\lambda_2-\lambda_1)}\int_{-\infty}^t (x(s)e^{-\lambda_2 s})^p e^{(p-1)\lambda_2 s}\, ds \leq \frac{1}{\Lambda\, (\lambda_2-\lambda_1)}\int_{-\infty}^t x(s)e^{-\lambda_2 s}e^{(p-1)\lambda_2 s}\, ds\, .
$$
Again Gronwall's inequality applied to the function $V(t)=\int_{-\infty}^t x(s)e^{-\lambda_2 s}e^{(p-1)\lambda_2 s}\, ds$ yields
$$
V(t_1)\leq V(t_0)\,  \exp \left( \frac{e^{(p-1)\lambda_2 t_1}}{\Lambda\, (\lambda_2-\lambda_1) (p-1)\lambda_2}\right)
$$
for all $t_0<t_1<T$. By letting $t_0\to -\infty$, since $V(t_0)\to 0$, we obtain $V(t_1)\leq 0$. But this  yields the contradiction $x(t)\equiv 0$.

 \hfill$\Box$

\subsection{Proof of  2.}

Let us consider now the case $p^\star\leq p <p^{\star \star}$, i.e. $\tau^- < \frac{2}{p-1}\leq \tau^+$.

By Proposition \ref{monotonicity}, we have that
$$
v(r)=r^{\tau^-}u(r)
$$
 is non decreasing. Hence, the limit $\lim_{r\to 0} u(r)r^{\tau^-}= \lim_{r\to 0} v(r) =v(0)$ exists and it is finite,  and the conclusion will follow by showing  that $v(0)>0$.

Since $\tau^-< {2\over p-1}$, we have  in particular that $\lim_{r\to 0} u(r)r^{\frac{2}{p-1}}=0$. By Proposition \ref{signeuprime}, this implies  that,  for $r>0$ small enough, one has $u^\prime \leq 0$ and  $u^{\prime \prime} \geq 0$, so that $u$ satisfies \eqref{equ}. Hence, by a direct computation, we deduce that $v$ is a nondecreasing solution of
\begin{equation}\label{eqv}
v''+\left( 1+2(\tau-\tau^-)\right) \frac{v'}{r}= \frac{v^p}{\Lambda r^{(p-1) \tau^-}} \quad \hbox{ for } 0<r\leq r_0\, ,
\end{equation}
 for some positive $r_0<1$.  
 Let us set $a\, := 1+2(\tau-\tau^-)>1$ and let us rewrite \eqref{eqv} as
 \begin{equation}\label{eqv1}
 \left( r^a v'(r)\right)' = \frac{1}{\Lambda} \frac{v^p}{r^{(p-1) \tau^--a}}>0\, .
 \end{equation}
 We emphasize that $a>1$ and $(p-1) \tau^-<2$.
 From \eqref{eqv1} it  follows in particular that the limit $\lim_{r\to 0} r^av'(r)=l\geq 0$ exists. If $l>0$, we get that $v'(r)\geq l r^{-a}$, which implies that $v$ is unbounded.
 Thus, $l=0$. Integrating \eqref{eqv1} and using that $v$ is bounded, we obtain
 $$
 v'(r)\leq c_1 r^{1-(p-1)\tau^-}
 $$
 for some $c_1>0$. Now, if, by contradiction, $v(0)=0$, we deduce by integration that
 $$
 v(r)\leq c_2 r^{2-(p-1)\tau^-}\, .
 $$
Using this estimate in \eqref{eqv1} and integrating again, we obtain a better upper bound on $v(r)$. Iterating the  procedure, we then conclude
$$
v(r)\leq c_k r^k\quad \hbox{ for every }  k>1\, .
$$
But this is a contradiction to Proposition \ref{bounded} as soon as $k\geq \tau^-$.
 
 \hfill$\Box$

\begin{rema}
{\rm The proof of statement 2 in Theorem \ref{behave} could be obtained by arguing as in the proof of statement 1 and by using also in this case the Emden-Fowler transformed $x(t)$. Indeed, Proposition \ref{monotonicity} gives that $x(t)\to 0$ as $t\to -\infty$ and, since the only positive root of the 
characteristic equation associated with the linearized equation is the eigenvalue $\lambda_2$ defined in \eqref{l1l2}, one can conclude that $x(t)\sim e^{\lambda_2 t}$ as $t\to -\infty$. 
We provided above an alternative elementary proof based on the properties of the solutions of equation \eqref{eqv}.}
\end{rema}

 \subsection{ Proof of 3.}                             
Let us now focus on  the case $p= p^{\star \star}$.                                    

We argue as in the proof of statement 2, and we consider the function
$$
v(r)=r^{\tau^-}u(r)\, ,
$$
which satisfies  now 
$$
v''+\left( 1+2(\tau-\tau^-)\right) \frac{v'}{r}= \frac{1}{\Lambda } \frac{v^p}{r^2} \quad \hbox{ for } 0<r\leq r_0
$$
 and, by Proposition \ref{p**},     
 $$
 v(0)=0\, .
 $$
 For $\bar K = \left[\Lambda \tau^- (\tau -\tau^-)\right]^{\tau^-/2}$ as in the statement, let us consider the function
 $$
 \underline{v}(r)= \frac{\bar K}{\left( -\log r +c_1\right)^{\frac{\tau^-}{2}}}\, ,
 $$
 with
 $$
 c_1= \left[ \left( \frac{\bar K}{v(r_0)}\right)^{\frac{2}{\tau^-}}+\log r_0\right]^+\, .
 $$
By a direct computation, we easily obtain
$$
\underline{v}''+\left( 1+2(\tau-\tau^-)\right) \frac{\underline{v}'}{r}\geq \frac{\tau^- (\tau-\tau^-) \bar K}{r^2 \left( -\log r +c_1\right)^{p \frac{\tau^-}{2}}}= \frac{1}{\Lambda} \frac{\underline{v}^p}{r^2}\, .
$$
Moreover, one has $v(0)=\underline{v}(0)=0$ and $v(r_0)\geq \underline{v}(r_0)$ by construction. This implies that 
$$
\min_{0\leq r\leq r_0} (v-\underline{v})(r)\geq 0\, .
$$
Indeed, if not, there exists an interior minimum point $\rho$ at which $v(\rho)<\underline{v}(\rho)$, $v'(\rho)=\underline{v}'(\rho)$ and $v''(\rho)\geq \underline{v}''(\rho)$. But this yields the contradiction
$$
\frac{1}{\Lambda} \frac{v^p(\rho)}{\rho^2}=v''(\rho) +\left( 1+2(\tau-\tau^-)\right) \frac{v'(\rho)}{\rho}
    \geq \underline{v}''(\rho) +\left( 1+2(\tau-\tau^-)\right) \frac{\underline{v}'(\rho)}{\rho}\geq \frac{1}{\Lambda} \frac{\underline{v}^p(\rho)}{\rho^2}\, .
    $$
Hence, we have obtained that
\begin{equation}\label{lower}
\underline{v}(r)= \frac{\bar K}{\left( -\log r +c_1\right)^{\frac{\tau^-}{2}}}\leq v(r)\quad \hbox{ for } 0\leq r\leq r_0\, .
\end{equation}
On the other hand, for any $b>\bar K$ and $c_2< -\log r_0$,  the function
$$
\overline{v}(r)=  \frac{b}{\left( -\log r -c_2\right)^{\frac{\tau^-}{2}}} 
$$
satisfies
$$
\begin{array}{ll}
\displaystyle \overline{v}''+\left( 1+2(\tau-\tau^-)\right) \frac{\overline{v}'}{r} & \displaystyle =\frac{b}{r^2} \left[ \frac{{\bar K}^{p-1}}{\Lambda (-\log r-c_2)^{\frac{\tau^-}{2}+1}}+ \frac{\tau^-(\tau^-+2)}{4 (-\log r-c_2)^{\frac{\tau^-}{2}+2}}\right]\\[3ex]
& \displaystyle  \leq \frac{b^p}{\Lambda r^2 (-\log r-c_2)^{\frac{\tau^-}{2}+1}}= \frac{\overline{v}^p}{\Lambda r^2}\quad \hbox{ for }0<r<r_0
\end{array}
$$
provided that
$$
b^{p-1}-{\bar K}^{p-1}\geq \frac{\Lambda \tau^-(\tau^-+2)}{4 (-\log r_0-c_2)}\, .
$$
We set
$$
c_2= -\log r_0 - \left( \frac{b}{v(r_0)}\right)^{\frac{2}{\tau^-}}\, ,
$$
and we fix $r_0$ so small that
 $$ \Lambda \tau^-( \tau^-+2) v(r_0)^{p-1} < b^{ p-1} ( b^{p-1}-\tilde K^{ p-1})\, .$$
With such a choice,  $\overline{v}$ is a super-solution for $0<r<r_0$, further satisfying $\overline{v}(0)=v(0)=0$ and $\overline{v}(r_0)=v(r_0)$. Arguing as above, we easily deduce that 
\begin{equation}\label{upper}
\overline{v}(r)= \frac{b}{\left( -\log r -c_2\right)^{\frac{\tau^-}{2}}}\geq v(r)\quad \hbox{ for } 0\leq r\leq r_0\, .
\end{equation}
Combining the bounds \eqref{lower} and \eqref{upper} and letting $r\to 0$, we obtain the conclusion by the arbitrariness of $b>\bar K$.

\hfill$\Box$  
                                  
\subsection{ Proof of 4.}
Let us finally consider the case $p>p^{\star \star}$. Arguing as in the proof of case 1, we distinguish two possible cases and we apply Proposition \ref{signeuprime}: either $\limsup_{r\to 0} u(r)r^{\frac{2}{p-1}}>0$, and then $u$ satisfies \eqref{2surp-1}, or    $\lim_{r\to 0} u(r)r^{\frac{2}{p-1}}=0$.  We claim that the second case cannot occur. Indeed, otherwise, we can consider the Emden--Fowler transformed $x(t)$ of $u$ as in the proof of case 1. Then, $x(t)$ is a solution of \eqref{eqx} satisfying $\lim_{t\to -\infty} x(t)=0$. However, under the current assumption on $p$, both the  eigenvalues $\lambda_1$ and $\lambda_2$ defined in \eqref{l1l2} are  negative, so that equation \eqref{eqx} does not posses  solutions vanishing as  $t\to -\infty$.

\hfill$\Box$


\begin{thebibliography}{99}
  
 \bibitem{BDLarxiv} I. Birindelli, F. Demengel, F. Leoni,  {\em Principal eigenvalues and eigenfunctions for fully nonlinear equations in punctured balls}, J. Math. Pures Appl. \textbf{186} (2024),  74--102. 

\bibitem{BGLP} I. Birindelli, G. Galise, F. Leoni, F. Pacella, \emph{Concentration and energy invariance for a class of fully nonlinear elliptic equations},  Calc. Var. Partial Differential Equations \textbf{57} (2018), no. 6, Art. 158, 22 pp.


\bibitem{CC} X. Cabr\'e, L. Caffarelli, Fully nonlinear elliptic equations, American Mathematical Society Colloquium Publications vol. 43, 1995. 
                          
 \bibitem{Ci} F.C. Cirstea,  \emph{A complete classification of the isolated singularities for nonlinear elliptic equations with inverse square potentials}, Memoirs of the American Mathematical Society \textbf{227}(1068) (2014), 1--97.
 
 \bibitem{CirsChau} F.C.  Cirstea, N. Chaudhuri, \emph{On trichotomy of positive singular solutions associated with the Hardy-Sobolev operator}, C. R. Acad. Sci. Paris, Ser. I \textbf{347} (3-4) (2009), 153--158.

 \bibitem{CDu1} F.C. Cirstea, Y. Du, \emph{Isolated singulariities for weighted quasilinear elliptic equations}, Journal of Functional Analysis \textbf{259} (1) (2010), 174--202.
 

\bibitem{CDu2} F.C.  Cirstea, Y. Du, \emph{ Asymptotic behavior of solutions of semilinear elliptic equations near an isolated singularity}, Journal of Functional Analysis  \textbf{250} (2) (2007),  317--346.


\bibitem{FQ} P. L. Felmer, A. Quaas, \emph{On critical exponents for the Pucci's extremal operators}, Ann. Inst. H. Poincar\'e Anal. Non Lin\'eaire  \textbf{20} (2003), no. 5, 843--865.

\bibitem{FS} W. Fleming, H. Soner,  Controlled Markov Processes and Viscosity Solutions, 
Springer-Verlag, New York, 1992.

\bibitem{F} R.H. Fowler, \emph{ Further studies on Emden's and similar differential equations},  Q. J. Math.  \textbf{2} (1931), 259--288.

\bibitem{Har} P. Hartman, Ordinary Differential Equations, 2nd ed., Birkh\"auser, Boston, 1982.

\bibitem{GIL} G. Galise, A. Iacopetti, F. Leoni,  \emph{Liouville-type results in exterior domains for radial solutions of fully nonlinear equations}, J. Differential Equations \textbf{269}, no. 6 (2020), 5034--506.

\bibitem{GILP} G. Galise, A. Iacopetti, F. Leoni,  F. Pacella, \emph{New concentration phenomena for a class of radial fully nonlinear equations}, Ann. Inst. H. Poincar\'e Anal. Non Lin\'eaire \textbf{37}, no. 5 (2020), 1109--1141.

\bibitem{L} P.L. Lions, \emph{Optimal control of diffusion processes and Hamilton-Jacobi-Bellman equations, part 2},  Comm. Partial Differential Equations \textbf{8} (1983), 1229--1276.

\bibitem{Pucci} C. Pucci, \emph{Operatori Ellittici Estremanti}, Ann. Mat. Pura Appl. \textbf{72} (1966), 141--170.

 \bibitem {WDu} L. Wei, Y. Du, \emph{Exact singular behavior of positive solutions to nonlinear elliptic equations with a Hardy potential}, J. Differential Equations  \textbf{262} (2017), 3864--3886.
\end{thebibliography}
 \end{document}